\documentclass[12pt, reqno]{amsart}
\usepackage{amsfonts}
\usepackage{amssymb}
\usepackage{amsmath}

\newtheorem{theorem}{Theorem}
\theoremstyle{plain}

\newtheorem{corollary}[theorem]{Corollary}

\newtheorem{lemma}[theorem]{Lemma}

\newtheorem{proposition}[theorem]{Proposition}

\begin{document}
\title{Identities involving values of Bernstein, $q$-Bernoulli and $q$-Euler polynomials}
\author{A. Bayad}
\address{Abdelmejid Bayad. D\'epartement de math\'ematiques \\
Universit\'e d'Evry Val d'Essonne, Bd. F. Mitterrand, 91025 Evry
Cedex, France  \\} \email{abayad@maths.univ-evry.fr}
\author{T. Kim}
\address{Taekyun Kim. Division of General Education-Mathematics \\
Kwangwoon University, Seoul 139-701, Republic of Korea  \\}
\email{tkkim@kw.ac.kr}

\maketitle

{\footnotesize {\bf Abstract} \hspace{1mm} {
In this paper we give some relations involving values of $q$-Bernoulli, $q$-Euler and Bernstein polynomials.  From these relations, we obtain some interesting identities on the $q$-Bernoulli, $q$-Euler and Bernstein polynomials.}

\section{Introduction and preliminaries }
 Throughout this paper, let $p$ be a fixed odd prime number. The symbols,  $\Bbb Z_p$, $\Bbb Q_p$ and $\Bbb C_p$ denote  the ring of $p$-adic integers, the field of $p$-adic numbers and the field of $p$-adic completion of the algebraic closure of $\Bbb Q_p$, respectively. When one talks of $q$-extension, $q$ is variously considered as an indeterminate, a complex number $q\in\Bbb C,$ or a $p$-adic $q\in\Bbb C_p$. Let $\Bbb N$ be the set of natural numbers and $\Bbb Z_+ =\Bbb N\cup \{0\}$.\\
Let $x=p^r\frac{s}{t}$ where $r\in\Bbb Q$ and $(s,t)=(p,s)=(p,t)=1.$ Then the $p$-adic absolute value is defined by $|x|_p=p^{-r}.$ If $q\in \Bbb C$,  we assume $|q|<1$, and if  $q\in \Bbb C_p$, we always assume $|1-q|_p <1$.\\
Let $C_{p^n}=\{\xi~|\xi^{p^n}=1\}$ be the cyclic group of order
$p^n$. Then,  the $p$-adic locally constant space is defined as $T_p=\lim_{n\to\infty}C_{p^n}$ ( see \cite{12}).
Let $UD(\Bbb Z_p)$ be the space of  uniformly differentiable functions  on $\Bbb Z_p$. For $f\in UD(\Bbb Z_p)$, the bosonic $p$-adic invariant integral on $\Bbb Z_p$ is defined by  
\begin{eqnarray}
I_1 (f)=\int_{\Bbb Z_p } f(x) d\mu (x)&=&\lim_{N \rightarrow \infty}  \sum_{x=0}^{p^N-1} f(x)\mu(x+p^N\Bbb Z_p)\\
&=&\lim_{N \rightarrow \infty}  \frac1{p^N}\sum_{x=0}^{p^N-1} f(x), \quad (\text{see \cite{9,12,13,14}}).\nonumber
\end{eqnarray}
From (1), we note that
\begin{eqnarray}
I_1(f_1)=\int_{\Bbb Z_p } f(x+1) d\mu(x)=\int_{\Bbb Z_p } f(x) d\mu(x)+f^{\prime}(0),
\end{eqnarray}
where $f_1 (x)=f(x+1)$, (see \cite{10,12}). \\
Let $f\in UD(\Bbb Z_p)$. Then the fermionic integral on $\Bbb Z_p$ is given by
\begin{eqnarray}
I_{-1}(f)=\int_{\Bbb Z_p } f(x) d\mu_{-1} (x)&=&
\lim_{N \rightarrow \infty}  \sum_{x=0}^{p^N-1} f(x)\mu_{-1}(x+p^N\Bbb Z_p)\\
&=&
\lim_{N \rightarrow \infty}  \sum_{x=0}^{p^N-1} f(x)(-1)^x, \quad (\text{see \cite{5,9,11}}).\nonumber
\end{eqnarray}
As known results, by (3), we get
\begin{eqnarray}
I_{-1}(f_1)=-I_{-1}(f)+2f(0), ~~(\textrm{ see \cite{6}}).
\end{eqnarray}
By using (2) and (4), we investigate some properties for the $q$-Bernoulli polynomials and the $q$-Euler polynomials.\\
Let $C[0, 1]$ be denote by the space of continuous functions on $[0, 1]$. For $f \in C[0, 1]$, Bernstein introduced the following well-known linear positive operator in the field of real numbers $\Bbb R$ :
\begin{eqnarray*}
\mathbb{B}_n(f | x)= \sum_{k=0}^{n} f(\frac{k}{n})\binom{n}{k} x^k
(1-x)^{n-k}=\sum_{k=0}^{n} f(\frac{k}{n})B_{k,n}(x),
\end{eqnarray*}
where
$\binom{n}{k}=\frac{n(n-1)\cdots(n-k+1)}{k!}=\frac{n!}{k!(n-k)!}$
(see \cite{1, 2, 5, 9, 10, 14}). Here, $\mathbb{B}_n (f|x)$ is called the
Bernstein operator of order $n$ for $f$. For $k,n \in \Bbb Z_+$, the Bernstein polynomials of degree $n$ are defined by
\begin{eqnarray}
B_{k,n}(x)= \binom{n}{k}x^k (1-x)^{n-k}, \quad \text{for} \ x\in [0,1].
\end{eqnarray}
For $x\in\Bbb Z_p$, the $p$-adic extension of Bernstein polynomials are given by 
\begin{eqnarray}
B_{k,n}(x)= \binom{n}{k}x^k (1-x)^{n-k}, \textrm{ where } k,n \in \Bbb Z_+,~~(\textrm{ see \cite{5,20,21}}.
\end{eqnarray}
The purpose of this paper is to give some relations for the $q$-Bernoulli, $q$-Euler and Bernstein polynomials.\\ From these relations, we obtain some interesteing identities on the $q$-Bernoulli, $q$-Euler and Bernstein polynomials.
\section{Some identities on the $q$-Bernoulli, $q$-Euler and Bernstein polynomials}
In this section we assume that $q\in T_p$. Let $f(x)=q^xe^{xt}.$ From (1) and (2), we have 
\begin{eqnarray}
\int_{\Bbb Z_p } f(x) d\mu (x)=\frac{t}{qe^t -1},~~(\textrm{ see \cite{12}}). 
\end{eqnarray}
The $q$-Bernoulli numbers are defined by
\begin{eqnarray}
\frac{t}{qe^t -1} =e^{B(q)t}=\sum_{n=0}^{\infty}
B_{n}(q)\frac{t^n}{n!}, ~(\textrm{ see \cite{12}}),
\end{eqnarray}
with usual convention about replacing $B_n(q)$ by $B_{n}(q).$ \\
By (7) and (8), we get Witt's formula for the $q$-Bernoulli numbers as follows:
\begin{eqnarray}
\int_{\Bbb Z_p} q^xx^n d\mu(x)=B_n(q),~n\in\Bbb Z_+,~(\textrm{ see \cite{12}}).
\end{eqnarray}
From (8), we can derive the following recurrence sequence:
\begin{eqnarray}
B_0(q)=0, \text{ and }  q(B(q)+1)^n-B_n(q)=\left\{\begin{array}{ll} 1 \ \ &\hbox{if}\ \ n=1,
\vspace{2mm}\\
0 \ \
&\hbox{if}\ \ n>1,
\end{array}\right.  
\end{eqnarray}
with usual convention about replacing $B^n(q)$ by $B_n(q)$.\\
Now, we define the $q$-Bernoulli polynomials as follows:
\begin{eqnarray}
\frac{t}{qe^t -1}e^{xt} =e^{B(x|q)t}=\sum_{n=0}^{\infty} B_{n}(x|q)\frac{t^n}{n!},
\end{eqnarray}
with usual convention about replacing $B^n(x|q)$ by $B_{n}(x|q).$ \\
By (2) and (11), we easily get
\begin{eqnarray}
\int_{\Bbb Z_p } q^ye^{(x+y)t} d\mu(y)=\frac{t}{qe^t-1}e^{xt}=\sum_{n=0}^{\infty} B_{n}(x|q)\frac{t^n}{n!}. 
\end{eqnarray}
Thus, we also obtain Witt's formula for the $q$-Bernoulli polynomials as follows:
\begin{eqnarray}
\int_{\Bbb Z_p } q^y(x+y)^n d\mu(y)= B_{n}(x|q),~\textrm{ for } n\in\Bbb Z_+. 
\end{eqnarray}
By (9) and (13), we easly get
\begin{eqnarray*}
B_n(x|q)&=&\sum_{l=0}^n\binom{n}{l}x^{n-l}B_l(q)\\
&=&\left(x+B(q)\right)^n,~\textrm{ for }n\in\Bbb Z_+.
\end{eqnarray*}
It is easy to show that
 \begin{eqnarray}
\frac{t}{qe^t -1}&=& \left(\frac{t}{q\left(1-q^{-1}\right)}\right) \left(\frac{1-q^{-1}}{e^t-q^{-1}}\right)\\
&=&\frac1{q\left(1-q^{-1}\right)}\sum_{n=0}^{\infty}H_n(q^{-1})(n+1)\frac{t^{n+1}}{(n+1)!},\notag
\end{eqnarray}
where $H_n(q^{-1})$ are the $n$-th Frobenius-Euler numbers ( see \cite{9}).\\
By (8), (10) and (15), we get 
\begin{eqnarray}
\frac{B_{n+1}(q)}{n+1}=\frac{H_n(q^{-1})}{q\left(1-q^{-1}\right)}=\frac{H_n(q^{-1})}{q-1} .
\end{eqnarray}
Therefore, by (15), we obtain the following proposition.
\begin{proposition}
For any $n\in\Bbb Z_+$, we have
\begin{eqnarray*}
\frac{B_{n+1}(q)}{n+1}=\frac{H_n(q^{-1})}{q\left(1-q^{-1}\right)}=\frac{H_n(q^{-1})}{q-1} ,
\end{eqnarray*}
where $H_n(q^{-1})$ are the $n$-th Frobenius-Euler numbers.
\end{proposition}
From (11) and (12), we can derive the following equation:
\begin{eqnarray}
\frac{qt}{qe^t -1}e^{(1-x)t}=\frac{-t}{q^{-1}e^{-t} -1}e^{-xt}=\sum_{n=0}^{\infty} B_{n}(x|q^{-1})(-1)^n\frac{t^n}{n!},
\end{eqnarray}
and
\begin{eqnarray}
\frac{qt}{qe^t -1}e^{(1-x)t}=qe^{B(1-x|q)t}=q\sum_{n=0}^{\infty} B_{n}(1-x|q)\frac{t^n}{n!}.
\end{eqnarray}
By comparing the coefficients on the both sides of (16) and (17), we obtain the following theorem.
\begin{theorem} Let $n\in\Bbb Z_+$. Then  we have
\begin{eqnarray*}
qB_{n}(1-x|q)=(-1)^nB_n(x|q^{-1}).
\end{eqnarray*}
\end{theorem}
 From (10) and (13), we have
\begin{eqnarray}
B_n(2|q)&=&\left(B(q)+1+1\right)^n=\sum_{l=0}^n\binom{n}{l}B_l(1|q)\\
&=&B_0(q)+\frac1{q}\sum_{l=1}^n\binom{n}{l}qB_l(1|q)=\frac1{q}\sum_{l=1}^n\binom{n}{l}B_l(q)\notag\\
&=&\frac1{q}B_n(1|q)=\frac1{q^2}qB_n(1|q)\notag=\frac1{q^2}B_n(q) \textrm{, while } n>1.\notag
\end{eqnarray}
Therefore, by (18), we obtain the following theorem.
\begin{theorem} For $n\in\Bbb Z_+$ with $n>1$, we have
\begin{eqnarray*}
q^2B_{n}(2|q)=B_n(q).
\end{eqnarray*}
\end{theorem}
By (9) and (13), we easily see that
\begin{eqnarray}
\int_{\Bbb Z_p} q^{-x}(1-x)^n d\mu(x)
&=& (-1)^n\int_{\Bbb Z_p} q^{-x}(x-1)^n d\mu(x)\\
&=&q\int_{\Bbb Z_p} (x+2)^nq^{x} d\mu(x)\notag\\
&=&\frac1{q}\int_{\Bbb Z_p} x^nq^{x} d\mu(x).\notag
\end{eqnarray}
Therefore, by (19), we obtain the following theorem.
\begin{theorem}
For $n \in \Bbb N$ with $n>1$, we have
\begin{eqnarray*}
\int_{\Bbb Z_p} q^{-x}(1-x)^n d\mu(x)=\frac1{q}\int_{\Bbb Z_p} x^nq^{x} d\mu(x).
\end{eqnarray*}
\end{theorem}
Let $f(x)=q^xe^{xt}.$  Then, from (3) and (4), we note that
  \begin{eqnarray} 
\int_{\Bbb Z_p} q^xe^{xt} d\mu_{-1}(x)=\frac{2}{qe^t +1},~(\textrm{ see \cite{5,17,18,19}}).
\end{eqnarray}
In \cite{5}, the $q$-Euler numbers are defined by 
\begin{eqnarray}
\frac{2}{qe^t +1} =e^{E(q)t}=\sum_{n=0}^{\infty} E_n(q)\frac{t^n}{n!},
\end{eqnarray}
with usual convention about replacing $E^n(q)$ by $E_n(q)$.\\
Let us define the $q$-Euler polynomials as follows:
\begin{eqnarray} 
\frac{2}{qe^t +1}e^{xt} =e^{E(x|q)t}=\sum_{n=0}^{\infty} E_n(x|q)\frac{t^n}{n!}.
\end{eqnarray}
From (3) and (4), we note that 
\begin{eqnarray} 
\int_{\Bbb Z_p } q^ye^{(x+y)t} d\mu_{-1} (y)=\frac{2}{qe^t +1}e^{xt}=\sum_{n=0}^{\infty} E_n(x|q)\frac{t^n}{n!}.
\end{eqnarray}
Thus, by (22) and (23), we get
\begin{eqnarray}
\int_{\Bbb Z_p} q^y(x+y)^n d\mu_{-1}(y)=E_n(x|q),~\textrm{ for } n\in\Bbb Z_+.
\end{eqnarray}
By (20), (21) and (24), we get
\begin{eqnarray}
E_n(x|q)=\sum_{l=0}^{n}\binom{n}{l}x^{n-l}E_{l,q}=\left(x+E(q)\right)^n,
\end{eqnarray}
with usual convention  $E(q)^n$ by $E_n(q).$\\
In \cite{22}, it is known that
\begin{eqnarray}
qE_n(2|q)=2+\frac1{q}E_n(q),
\end{eqnarray}
and
\begin{eqnarray*}
\int_{\Bbb Z_p }q^{-x} (1-x)^n d\mu_{-1} (x)=
2+\frac1{q}\int_{\Bbb Z_p }x^n q^x d\mu_{-1} (x),~( \textrm{ see \cite{7}}).
\end{eqnarray*}
For $x\in\Bbb Z_p$, by (6), we have
\begin{eqnarray}
B_n(1-x|q)&=&\sum_{l=0}^n\binom{n}{l}(1-x)^{n-l}B_l(q)\\
&=& \sum_{l=0}^n\binom{n}{l}(1-x)^{n-l}x^lB_l(q)x^{-l}\notag \\
&=& \sum_{l=0}^n B_{l,n}(x)B_l(q)x^{-l},\notag
\end{eqnarray}
where $B_{l,n}(x)$ are Bernstein polynomials of degree $n$.\\
Therefore, by (27), we obtain the following proposition.
\begin{proposition}For $n\in\Bbb Z_+,$ we have
\begin{eqnarray*}
B_n(1-x|q)=\sum_{l=0}^n B_{l,n}(x)B_l(q)x^{-l},
\end{eqnarray*}
\end{proposition}
where $B_{l,n}(x)$ are Bernstein polynomials of degree $n$.\\
From (7) and (20), we can derive the following equation:
\begin{eqnarray} 
\int_{\Bbb Z_p }\int_{\Bbb Z_p } q^{x+y}e^{(x+y)t}d\mu_{-1} (x) d\mu (y)=\frac{2t}{q^2e^{2t}-1}.
\end{eqnarray}
By (8), we get 
\begin{eqnarray} 
\frac{2t}{q^2e^{2t}-1}=e^{B(q^2)2t}=\sum_{n=0}^{\infty} 2^nB_n(q^2)\frac{t^n}{n!}.
\end{eqnarray}
By comparing the coefficients on the both sides of (28) and (29), we obtain the following theorem.
\begin{theorem}For $n\in\Bbb Z_{+}$, we have 
\begin{eqnarray*}
\frac1{2^n} \int_{\Bbb Z_p }\int_{\Bbb Z_p } q^{x+y}(x+y)^{n+1} d\mu_{-1} (x) d\mu (y)=B_n(q^2).
\end{eqnarray*}
\end{theorem}
By using (9) and (24), we obtain the following corollary.
\begin{corollary}
For $n \in \Bbb Z_+$, we have
\begin{eqnarray*}
\frac1{2^n} \displaystyle\sum_{l=0}^{n} \binom{n}{l} E_l(q)B_{n-l}(q)=B_n(q^2). 
\end{eqnarray*}
\end{corollary}
From the definition of Bernstein polynomials, we note that
\begin{eqnarray*}
B_{k,n}(x)=B_{n-k,n}(1-x),~\textrm{ for } n, k \in \Bbb Z_+.
\end{eqnarray*}
Thus, we have
\begin{eqnarray*}
B_{k,n}(\frac12 )=B_{n-k,n}( \frac12 )=\binom{n}{k} \left(\frac12\right)^{n-k} \left(\frac12 \right)^{k} =\left(\frac12\right)^{n}\binom{n}{k}.
\end{eqnarray*}
Therefore, we obtain the following lemmma.
\begin{lemma}Let $n, k \in \Bbb Z_+$. Then we have 
\begin{eqnarray*}
B_{k,n}(\frac12 )=\left(\frac12\right)^{n}\binom{n}{k}.
\end{eqnarray*}
\end{lemma}
By Lemma 8 and Corollary 7, we obtain the following corollary.
\begin{corollary}
For $n \in \Bbb Z_+$, we have
\begin{eqnarray*}
B_n(q^2)&=& \frac1{2^n} \displaystyle\sum_{l=0}^{n} \binom{n}{l} E_l(q)B_{n-l}(q)\\
&=& \displaystyle\sum_{l=0}^{n} B_{l,n}(\frac12 ) E_l(q)B_{n-l}(q).
\end{eqnarray*}
\end{corollary}
For the right side of (28), we have
\begin{eqnarray} 
\frac{2t}{q^2e^{2t}-1}&=&-2t\sum_{m=0}^{\infty}q^{2m}e^{2mt}=
-2t\sum_{n=0}^{\infty}\left(2^n \sum_{m=0}^{\infty}q^{2m}m^n\right)\frac{t^n}{n!}\\
&=& \sum_{n=0}^{\infty}\left(-2^{n+1}(n+1) \sum_{m=0}^{\infty}q^{2m}m^n\right)\frac{t^{n+1}}{(n+1)!}.\notag
\end{eqnarray}
In (29), we see that $B_0(q^2)=0.$ By comparing coefficients on the both sides of (29) and (30), we obtain the following theorem.
\begin{theorem} For $n\in\Bbb Z_+$, we have 
\begin{eqnarray*}
-\frac{1}{2^{(n+1)}(n+1)}
\int_{\Bbb Z_p }\int_{\Bbb Z_p } q^{x+y}(x+y)^nd\mu_{-1} (x) d\mu (y)= \sum_{m=0}^{\infty}q^{2m}m^n.
\end{eqnarray*}
\end{theorem}
By (9), (24) and binomial theorem, we obtain the following corollary
\begin{corollary}
For $n \in \Bbb Z_+$, we have
\begin{eqnarray*}
\sum_{m=0}^{\infty}q^{2m}m^n= -\frac1{2^{(n+1)}(n+1)} \displaystyle\sum_{l=0}^{n+1} \binom{n+1}{l} E_l(q)B_{n+1-l}(q).
\end{eqnarray*}
\end{corollary}
By Lemma 7, we obtain the following corollary.
\begin{corollary}
For $n \in \Bbb Z_+$, we have
\begin{eqnarray*}
\sum_{m=0}^{\infty}q^{2m}m^n= -\frac1{n+1} \displaystyle\sum_{l=0}^{n+1} B_{n+1,l}(\frac12 )E_l(q)B_{n+1-l}(q).
\end{eqnarray*}
\end{corollary}
It seems to be important to study double $p$-adic integral representation of bosonic and fermionic on the Bernstein polynomials associated with $q$-Bernoulli and $q$-Euler polynomials. Theorem 10 is useful to study those integral representation on $\Bbb Z_p$ related to Bernstein polynomials. Now, we take bosonic $p$-adic invariant integral on $\Bbb Z_p$ for one Bernstein polynomials in (6)as follows:
\begin{eqnarray}
\int_{\Bbb Z_p }B_{k,n}(x) q^{x} d\mu (x)&=&\int_{\Bbb Z_p }
\binom{n}{k}(1-x)^{n-k}x^kq^xd\mu (x)\\
&=&\binom{n}{k} \sum_{l=0}^{n-k}\binom{n-k}{l}(-1)^l\int_{\Bbb Z_p }x^{k+l}q^{x} d\mu (x)\notag\\
&=&\binom{n}{k} \sum_{l=0}^{n-k}\binom{n-k}{l}(-1)^lB_{k+l}(q).\notag
\end{eqnarray}
From the reflection symmetry of Bernstein polynomials, we have
\begin{eqnarray}
B_{k,n}(x)=B_{n-k,n}(1-x),~\textrm{ where } n, k\in\Bbb Z_{+}.
\end{eqnarray}
Let $n, k\in\Bbb N$ with $n>k+1$. Then, by (32), we get
\begin{eqnarray}
\int_{\Bbb Z_p }B_{k,n}(x) q^{x} d\mu (x)&=&\int_{\Bbb Z_p }B_{n-k,n}(1-x) q^{x} d\mu (x)\\
&=&\binom{n}{k}\sum_{l=0}^{k}\binom{k}{l}(-1)^{k-l} \int_{\Bbb Z_p }(1-x)^{n-l}q^xd\mu (x)\notag\\
&=&q\binom{n}{k}\sum_{l=0}^{k}\binom{k}{l}(-1)^{k-l} \int_{\Bbb Z_p }q^{-x}x^{n-l}d\mu (x)\notag\\
&=&q\binom{n}{k}\sum_{l=0}^{k}\binom{k}{l}(-1)^{k-l} B_{n-l}(q^{-1}).\notag
\end{eqnarray}
Therefore, we obtain the following theorem.
\begin{theorem} Let $n, k \in \Bbb Z_+$ with $n>k+1$. Then we have
\begin{eqnarray*}
\int_{\Bbb Z_p }B_{k,n}(x) q^{1-x} d\mu (x)= 2^k\binom{n}{k}\sum_{l=0}^k B_{l,k}(\frac12 )(-1)^{k-l}B_{n-l}(q).
\end{eqnarray*}
\end{theorem}
By (31), we obtain the following corollary
\begin{corollary}
For $n, k \in \Bbb Z_+$ with $n>k+1$, we have
\begin{eqnarray*}
2^{2k}\binom{n}{k} B_{l,k}(\frac12 )(-1)^{k-l}B_{n-l}(q)= 2^nq
\sum_{l=0}^{n-k}B_{l,n-k}(\frac12 )(-1)^{l}B_{k+l}(q).
\end{eqnarray*}
\end{corollary}
Let $m,n,k\in\Bbb Z_+$, with $m+n>k+1$. Then we have
\begin{eqnarray}
\int_{\Bbb Z_p }B_{k,n}(x)B_{k,m}(x) q^{-x} d\mu (x)
&=&\binom{n}{k}\binom{m}{k}\sum_{l=0}^{2k} \binom{2k}{l}(-1)^{l+2k} \int_{\Bbb Z_p }q^{-x}(1-x)^{n+m-l}d\mu (x)\notag\\
&=&\binom{n}{k}\binom{m}{k}\sum_{l=0}^{2k} \binom{2k}{l}(-1)^{l+2k} q\int_{\Bbb Z_p }q^{x}(x+2)^{n+m-l}d\mu (x)\notag\\
&=&\binom{n}{k}\binom{m}{k}\sum_{l=0}^{2k} \binom{2k}{l}(-1)^{l+2k} \frac1{q}B_{n+m-l}(q).
\end{eqnarray}
Therefore, by (34), we obtain the following theorem.
\begin{theorem} Let $m,n,k\in\Bbb Z_+$ with $m+n>k+1$. Then we have
\begin{eqnarray*}
\int_{\Bbb Z_p }B_{k,n}(x)B_{k,m}(x) q^{1-x} d\mu (x)=\binom{n}{k}\binom{m}{k}2^{2k}\sum_{l=0}^{2k} B_{l,2k}(\frac12 )(-1)^{l+2k} B_{n+m-l}(q).
\end{eqnarray*}
\end{theorem}
By binomial theorem, we easily get 
\begin{eqnarray}
\int_{\Bbb Z_p }B_{k,n}(x)B_{k,m}(x) q^{1-x} d\mu (x)
&=&\binom{n}{k}\binom{m}{k}\sum_{l=0}^{n+m-2k} \binom{n+m-2k}{l}(-1)^l\int_{\Bbb Z_p }x^{l+2k} q^{1-x} d\mu (x) )\notag\\
&=&q\binom{n}{k}\binom{m}{k}\sum_{l=0}^{n+m-2k} \binom{n+m-2k}{l}(-1)^l B_{l+2k}(\frac1q ).
\end{eqnarray}
By (35), we obtain the following corollary.
\begin{corollary}
Let $m,n, k \in \Bbb Z_+$ with $n+m>2k+1$. Then we have
\begin{eqnarray*}
2^{4k}\sum_{l=0}^{2k} (-1)^{l+2k} B_{l,2k}(\frac12)B_{n+m-l}(q )= 2^{n+m}q\sum_{l=0}^{n+m-2k} (-1)^l B_{l,n+m-2k}(\frac12)B_{l+2k}(\frac1q ).
\end{eqnarray*}
\end{corollary}
For $s \in\Bbb N$, let $n_1, n_2, \ldots, n_s, k \in \Bbb Z_+$ with $n_1+n_2+\cdots+ n_s>sk+1.$
\begin{eqnarray}
&&\int_{\Bbb Z_p} 
B_{k, n_1}(x)\cdots   B_{k, n_s}(x)  q^{-x}d\mu(x)
=
\binom{n_1}{k}\cdots \binom{n_s}{k} \int_{\Bbb  Z_p} x^{sk}(1-x)^{n_1+...+n_s}q^{-x}d\mu(x)\notag\\
=&&\binom{n_1}{k}\cdots \binom{n_s}{k} \sum_{l=0}^{sk}(-1)^{sk+l}\int_{\Bbb  Z_p} (1-x)^{n_1+...+n_s-l}q^{-x}d\mu(x)\notag\\
=&&q\binom{n_1}{k}\cdots \binom{n_s}{k} \sum_{l=0}^{sk}\binom{sk}{l}(-1)^{sk+l}\int_{\Bbb  Z_p} (x+2)^{n_1+...+n_s-l}q^{x}d\mu(x)\notag\\
=&&\frac1q \binom{n_1}{k}\cdots \binom{n_s}{k} \sum_{l=0}^{sk} \binom{sk}{l}(-1)^{sk+l}B_{n_1+\cdots+n_s-l}(q).
\end{eqnarray}
Therefore, by (36), we obtain the following theorem.
\begin{theorem} For $s \in\Bbb N$, let $n_1, n_2, \ldots, n_s, k \in
  \Bbb Z_+$ with $n_1+n_2+\cdots+ n_s>sk+1.$ Then we have
\begin{eqnarray*}
\int_{\Bbb Z_p}B_{k, n_1}(x)\cdots   B_{k, n_s}(x)  q^{1-x}d\mu(x)=
\binom{n_1}{k}\cdots \binom{n_s}{k} 2^{sk}\sum_{l=0}^{sk} (-1)^{sk+l}B_{l,sk}(\frac12 )B_{n_1+\cdots+n_s-l}(q).
\end{eqnarray*}
\end{theorem}
From binomial theorem, we can easily derive the following equation:
\begin{eqnarray}
&&\int_{\Bbb Z_p} B_{k, n_1}(x)\cdots   B_{k, n_s}(x)  q^{-x}d\mu(x)\\
&&=\binom{n_1}{k}\cdots \binom{n_s}{k} \sum_{l=0}^{n_1+\cdots+n_s-sk}\binom{n_1+\cdots+n_s-sk}{l}(-1)^l \int_{\Bbb  Z_p} x^{l+sk}q^{-x}d\mu(x)\notag\\
&&=\binom{n_1}{k}\cdots \binom{n_s}{k} \sum_{l=0}^{n_1+\cdots+n_s-sk}\binom{n_1+\cdots+n_s-sk}{l}(-1)^lB_{l+sk}(q^{-1}).\notag
\end{eqnarray}
By (37), we obtain the following corollary.
\begin{corollary} Let  $s \in\Bbb N$ and  $n_1, n_2, \ldots, n_s, k \in \Bbb Z_+$ with $n_1+n_2+\cdots+ n_s>sk+1.$ Then we have
\begin{eqnarray*}
2^{n_1+\cdots+n_s-2sk}q\sum_{l=0}^{n_1+\cdots+n_s-sk}(-1)^lB_{l,n_1+\cdots+n_s-sk}(\frac12)B_{l+sk}(q^{-1})\\
= \sum_{l=0}^{sk}B_{l,sk}(\frac12)(-1)^{l+sk}B_{n_1+\cdots+n_s-l}(q).
\end{eqnarray*}
\end{corollary}
\section{Further remarks}
In this ection, we assume that $q\in\Bbb C$ with $|q|<1.$ The $q$-Euler polynomials and $q$-Bernoulli polynomials are defined by the generating functions as follows:
 \begin{eqnarray*}
F_q^{E}(t,x)=\frac{2}{qe^t +1}e^{xt}=\sum_{n=0}^{\infty} E_{n}(q|x)\frac{t^n}{n!},~|t+logq|<\pi ,
\end{eqnarray*}
and
\begin{eqnarray}
F_q^{B}(t,x)=\frac{t}{qe^t -1}e^{xt}=\sum_{n=0}^{\infty} B_{n}(q|x)\frac{t^n}{n!},~|t+logq|<2\pi .
\end{eqnarray}
In the special case $x=0$, $E(0|q)=E_n(q)$ are called the $n$-th $q$-Euler numbers and $B_n(0|q)=B_n(q)$ are also called the $n$-th $q$-Bernoulli numbers. \\
As usual convention, let us define $F_q^{B}(t,0)=F_q^{B}(t)$ and $F_q^{E}(t,0)=F_q^{E}(t)$. \\ For $s\in\Bbb C$, we have 
\begin{eqnarray}
\frac{q}{\Gamma(s)}\int_{0}^{\infty}F_q^{B}(-t,1)t^{s-2}dt&=&\frac{q}{\Gamma(s)}\int_{0}^{\infty}\frac{e^{-t}}{1-qe^{-t}}t^{s-1}dt\\
&=&\sum_{m=0}^{\infty}q^{m+1}\frac{1}{\Gamma(s)}\int_{0}^{\infty}e^{-(m+1)t}t^{s-1}dt \notag\\
&=&\sum_{m=0}^{\infty}\frac{q^{m+1}}{(m+1)^s}=\sum_{m=1}^{\infty}\frac{q^m}{m^s}.\notag
\end{eqnarray}
From (39), we define$q$-zeta function as follows: 
\begin{eqnarray*}
\zeta_q(s)=\sum_{m=1}^{\infty}\frac{q^m}{m^s}, \textrm{ for } \Re(s)>1. 
\end{eqnarray*}
Note that $\zeta_q(s)$ has meromorphic continuation to the whole complex $s$-plane with a simple pole at $s=1$.\\
By (38), (39) and elementary complex integral, we get
\begin{eqnarray}
\zeta_q(1-n)=(-1)^n\frac{qB_n(1|q)}{n}=\left\{\begin{array}{ll} -(1+B_1(q))  \ \ &\hbox{ if }\ \ n=1,
\vspace{2mm}\\
(-1)^n\frac{B_n(q)}{n} \ \
&\hbox{if}\ \ n>1.
\end{array}\right.  
\end{eqnarray}
From (38), we note that
\begin{eqnarray*}
\left(\frac{2}{qe^t +1}\right) \left(\frac{t}{qe^t -1}\right) =\frac{2t}{q^2e^{2t} -1}=\sum_{n=0}^{\infty}2^n B_{n}(q^2)\frac{t^n}{n!},
\end{eqnarray*}
and
\begin{eqnarray*}
\left(\frac{2}{qe^t +1}\right) \left(\frac{t}{qe^t -1}\right) =e^{\left(B(q)+E(q)\right)t}=\sum_{n=0}^{\infty}\left(B(q)+E(q)\right)^n\frac{t^n}{n!}.
\end{eqnarray*}
Thus, we have
\begin{eqnarray}
2^nB_n(q^2)=(B(q)+E(q))^n=\sum_{l=0}^n\binom{n}{l}B_l(q)E_{n-l}(q),
\end{eqnarray}
with usual convention about replacing $B^n(q)$ by $B_n(q).$\\
By (42), we get 
\begin{eqnarray}
B_n(q^2)=\frac1{2^n}\sum_{l=0}^n\binom{n}{l}B_l(q)E_{n-l}(q), (\textrm{ cf.\cite{8}}).
\end{eqnarray}
From (40) and (41), we can derive the following equation:
\begin{eqnarray*}
\zeta_{q^2}(-n)&=&(-1)^{n+1}\frac{q^2B_{n+1}(q^2)}{n+1}=-(-1)^{n}\frac{B_{n+1}(q^2)}{n+1}\\
&=&-(-1)^{n}\frac1{2^{n+1}(n+1)}\sum_{l=0}^{n+1}\binom{n+1}{l}B_l(q^2)E_{n+1-l}(q^2).
\end{eqnarray*}

\end{document}